\magnification=\magstep1
\input amstex
\documentstyle{amsppt}
\TagsOnRight
\hsize=5in                                                  
\vsize=7.8in
\strut
\vskip5truemm
\font\small=cmr8
\font\itsmall=cmti8

\def\smallarea#1{\par\begingroup\baselineskip=10pt#1\endgroup\par}
\def\abstract#1{\begingroup\leftskip=5mm \rightskip 5mm
\baselineskip=10pt\par\small#1\par\endgroup}

\def\refer#1#2{\par\begingroup\baselineskip=11pt\keyindent\leftskip=8.25mm\rightskip=0mm
 \strut\llap{#1\kern 1em}{#2\hfill}\par\endgroup}
              
\nopagenumbers
\centerline{\bf COMBINATORIAL INVARIANTS COMPUTING}
\centerline{\bf THE RAY-SINGER ANALYTIC TORSION}
\vskip 25.3pt

\centerline{MICHAEL S. FARBER
\footnotemark"$^1$"}
\footnotetext"$^{1}$"{The research was supported
by grant No.~449/94-1 from the Israel Academy of Sciences and Humanities}
\smallarea{%
\centerline {\itsmall School of Mathematical Sciences,}
\centerline {\itsmall Tel-Aviv University,}
\centerline {\itsmall Tel-Aviv 69978, Israel}
\centerline {\itsmall e-mail farber $\@$math.tau.ac.il}
\centerline {\itsmall fax 972-3-6407543}}

 

\input amstex
\define\C{{\Bbb C}}
\define\R{{\Bbb R}}

\define\Hom{\operatorname{Hom}}


\redefine\D{\Cal D}

\def\<{\langle}
\def\>{\rangle}

\define\pd#1#2{\dfrac{\partial#1}{\partial#2}}

\documentstyle{amsppt}   
\vskip 1cm

Let $K$ denote a closed odd-dimensional smooth manifold
and let $E$ be a flat vector bundle over $K$. In this situation the 
construction of Ray and Singer \cite{RS} gives a metric on the determinant 
line of the cohomology $\det H^\ast(M;E)$ which is a {\it smooth invariant 
of the manifold $M$ and the 
flat bundle} $E$. (Note that if the dimension of $K$ is even then the 
Ray-Singer metric depends on the choice of a Riemannian metric on $K$ 
and of a Hermitian metric on $E$).
The famous theorem which was 
proved by J. Cheeger \cite{C} and W.M\"uller \cite{Mu}, states that 
assuming that the flat vector bundle $E$ is {\it unitary} (i.e. $E$ 
admits a flat Hermitian metric), the Ray-Singer metric 
coincides with the Reidemeister metric which is defined using finite
dimensional linear algebra by the
combinatorial structure of $K$. 

The construction of the Reidemeister metric works also 
under a weeker assumption that the flat bundle $E$ is {\it unimodular}, i.e.
the line bundle $\det(E)$ admits a flat metric. In a recent paper
W.M\"uller \cite{Mu1} proved that in this case the Ray-Singer metric 
again coincides with the Reidemeister metric.

Without the assumption that the flat bundle $E$ is unimodular,
the standard construction of Reidemeister metric is ambiguous
(it depends on different choices made in the process of its construction).
In this general situation J.-M. Bismut and W. Zhang \cite{BZ}
computed the deviation of the 
Ray-Singer metric from so called Milnor metric, cf. \cite{BZ}; 
the result
of their computation is given in a form of an integral of a Chern-Simons 
current which compensates the ambiguity of the Milnor metric and 
depends on the Riemannian metric on $K$ and on the metric on $E$.

In this paper I show that for any piecewise-linear closed orientable
manifold $K$ of odd dimension there exits an invariantly defined metric
on the determinant line of cohomology $\det(H^\ast(K;E))$, 
where $E$ is an arbitrary flat bundle over $K$. The construction of this
metric is purely combinatorial. I call this metric Poincar\'e - Reidemeister
metric since it is defined by combining the standard Reidemeister's
construction with the Poincar\'e duality. 

The idea to use the Poincar\'e duality in order to construct an invariantly
defined metric on the determinant line was
prompted by the theorem 4.1 of D.Burghelea, L.Friedlander, and 
T.Kappeler \cite{BFK}, expressing the analytic torsion through some
torsion invariants (depending on the Riemannian metric) associated with
a Morse function.

The main properties of the 
Poincar\'e - Reidemeister metric consist in the following: 

(a) the Poincar\'e-Reidemeister metric can be computed starting from any
polyhedral cell decomposition of the manifold and using purely combinatorial
terms, cf. 4.3.

(b) the Poincar\'e - Reidemeister metric coincides with the 
Reidemeister metric when the latter is correctly defined; 
 
(c) The construction of Ray and Singer, which uses zeta-function regularized
determinants of Laplacians, produces the metric on the determinant
of cohomology, which coincides (via the De Rham isomorphism) with the 
Poincar\'e-Reidemeister metric; this is the main result of the paper,
formulated as Theorem 6.2;

(d) The Poincar\'e-Reidemeister metric behaves well with respect to
natural correspondences between determinant lines which are discussed
in \S 1, cf. Prop. 4.8. 

Another interesting observation made here consists in finding that
Ray-Singer metrics on some relative determinant lines can be computed
combinatorially (including the even-dimensional case) in terms of the metrics
determined by correspondences, cf. Theorems 5.1, 5.2 and 5.3.

The work technically is based on the fundamental theorem (0.2) of 
J.-M.Bismut and W.Zhang 
\cite{BZ} and uses actually a very special corollary of it which is  
proved by cancelling all complicated terms in the Bismut-Zhang theorem
(cf. Theorem 5.1 below). Note, that even weeker theorem 5.3 allows to 
{\it identify completely the Ray-Singer norm (cf. proof of Theorem 6.2) in 
terms of the Poincar\'e-Reidemeister norm in the odd-dimensional 
case}. Most of the paper consists in careful elementary analysis of the 
foundations: we tried to study separately different duality relations 
which exist among the determinant lines,
duality between homology and cohomology and the Poincar\'e duality,
and also correspondences induced by isomorphisms of the volume bundles.

\subheading{Notations} $K$ will denote a finite polyhedron given by its 
polyhedral cell decomposition (cf. \cite{RS1}).
We will consider flat vector bundles $E$ over $K$;
this we understand in the following way: the structure group of $E$ has been
reduced to a descrete group. We will also identify flat vector bundles
with the locally constant sheaves of their flat sections.

Given a polyhedral cell decomposition $\tau$ of $K$, one constructs 
the chain complex
$C_\ast(M,\tau,E)$ as follows: the basis of $C_q(M,\tau,E)$ form
pairs $(\Delta^q;s)$, where $\Delta^q$ is an oriented $q$-dimensional
cell of $\tau$ and $s$ is a flat section of $E$ over $\Delta^q$.
The boundary operator 
$$\partial: C_q(M,\tau,E)\to C_{q-1}(M,\tau,E)$$
is defined by the usual formula
$$\partial(\Delta^q;s)=\sum \epsilon[\Delta^q,\Delta^{q-1}](\Delta^{q-1},
s^\prime),$$
where $\epsilon[\Delta^q,\Delta^{q-1}]$ is the sign $\pm 1$ determined by the
orientations and $s^\prime$ is the restriction of $s$ on $\Delta^{q-1}$ and 
the sum runs over all $(q-1)$-dimensional cells
$\Delta^{q-1}\subset\overline{\Delta^q}$.

We will denote by 
$$T: \det(C_\ast(M,\tau,E))\to \det(H_\ast(M;E))$$
the canonical isomorphism between the determinant lines, cf. \cite{BGS}. 
Sometimes we will
write $T_\tau$ instead of $T$ in order to emphasize dependance on the
polyhedral cell decomposition $\tau$.

\subheading{Warning} In this paper, like in most of the literature, we
neglect signs of isomorphism between determinant lines. To be more precise,
we will work in the category whose objects are one-dimensional vector
spaces (over $\R$ or over $\C$) and whose morphisms are linear maps,
such that $f$ and $-f$ are considered as representing the same morphism.

\heading 1. Correspondences between determinant lines\endheading

\subheading{1.1} Suppose that we are given two flat 
vector bundles $E$ and $F$ (real or complex)
over a finite polyhedron $K$ and an 
isomorphism of the determinant flat line bundles 
$\phi:\det(E)\to \det(F)$. 
We will show in this section
that these data determine canonically an isomorphism between the 
corresponding determinant lines
$$\hat\phi : \det H_\ast(K;E)\to \det H_\ast(K;F).$$
This map $\hat\phi$ will be called {\it correspondence 
between the determinant lines} determined by $\phi$.

Let $\tau$ be a polyhedral cell decomposition of $K$.
According to the definition above, 
the vector space $C_q(M,\tau,E)$ is the direct sum of the spaces of
flat sections $\Gamma(\Delta;E)$, where $\Delta$ runs over all oriented
$q$-dimensional cells of $K$. Thus one can identify the determinant line
$\det C_q(K,\tau,E)$ with the tensor product of determinant lines
$$\otimes_\Delta \det\Gamma(\Delta;E).$$
On the other hand, for any cell $\Delta$ 
we have the canonical isomorphism
$\det \Gamma(\Delta;E)=\Gamma(\Delta;\det(E)).$
Therefore we obtain the canonical isomorphisms
$$\det(C_\ast(K,\tau,E))=\otimes_q \det(C_q(K,\tau,E))^{(-1)^q}=
\otimes \Gamma(\Delta,\det(E))^{\epsilon(\Delta)},\tag1$$
where $\epsilon(\Delta)$ denotes $(-1)^{\dim(\Delta)}$.
We may consider the similar canonical decomposition for the flat bundle
$F$ as well and then one defines the map
$$\tilde\phi: \det(C_\ast(K,\tau,E))\to \det(C_\ast(K,\tau,F))\tag2$$
as the tensor product the maps 
$\Gamma(\Delta,\det(E))\to\Gamma(\Delta,\det(F))$, 
induced by the map $\phi$ if the dimension $\dim(\Delta)$ is even, and 
induced by the
contragradient map $\overline\phi:(\det(E))^\ast\to(\det(F))^\ast$
(inverse to the adjoint)
if the dimension $\dim(\Delta)$ is odd. Finally we define the map 
$\hat\phi$ which appeares in the commutative diagram
$$
\CD
\det(C_\ast(K,\tau,E))@>{\tilde\phi}>> \det(C_\ast(K,\tau,F))\\
@V{T_\tau}VV  @V{T_\tau}VV\\
\det(H_\ast(K,E))@>{\hat\phi}>> \det(H_\ast(K,F)) \\
\endCD
\tag3
$$
Recall that the vertical maps $T_\tau$ are the canonical maps, cf.\cite{BGS}.
The above definition of $\hat\phi$ is justified by combinatorial invariance
property, similar to combinatorial invariance of the Reidemeister torsion
\cite{M}:

\proclaim{1.2. Proposition (Combinatorial Invariance)} 
The map $\hat\phi$ does not depend on the polyhedral cell decomposition 
$\tau$ of $K$. More precisely, the above construction gives
the same map $\hat\phi$ between the determinant lines for any pair of
polyhedral cell decompositions of $K$ having a common subdivision.
\endproclaim
\demo{Proof} Let $\tau^\prime$ be a subdivision of the polyhedral cell 
decomposition $\tau$.
We have the natural inclusion of the chain complexes
$$i: C_\ast(K,\tau,E)\to C_\ast(K,\tau^\prime,E);$$
let $D_\ast(\tau^\prime,\tau,E)$ denote the factor complex. 
It is acyclic and so there is correctly defined a volume 
$$\alpha(\tau^\prime,\tau,E)\in\det(D_\ast(\tau^\prime,\tau,E)).\tag4$$

The spaces of chains of the complex 
$D_\ast=D_\ast(\tau^\prime,\tau,E)$ are the factor-spaces $D_q$
which appear in the
exact sequence
$$0\to\oplus_\Delta\Gamma(\Delta,E)@>r>>
\oplus_{\Delta^\prime}\Gamma(\Delta^\prime,E)\to D_q\to 0;$$
here $\Delta$ runs over all cells of $\tau$ which are
not cells of $\tau^\prime$ and $\Delta^\prime$ runs over all 
cells of $\tau^\prime$ which are not cells of $\tau$. The map $r$ is given
by restriction of flat sections. Therefore we obtain
$$\det(D_\ast(\tau^\prime,\tau,E))=
\otimes_{\Delta^\prime}
\Gamma(\Delta^\prime,\det(E))^{\epsilon(\Delta^\prime)}
\otimes_{\Delta} 
\Gamma(\Delta,\det(E))^{\epsilon(\Delta)+1}.\tag5$$
In the last formula again
$\Delta$ runs over all cells of $\tau$ which are
not cells of $\tau^\prime$ and $\Delta^\prime$ runs over all 
cells of $\tau^\prime$ which are not cells of $\tau$. 

Now suppose that we are given two flat bundles $E$ and $F$ over the complex
$K$ and a bundle
isomorphism $\phi:\det(E)\to\det(F)$. We may construct the corresponding
factor-complexes $D_\ast(\tau^\prime,\tau,E)$ and 
$D_\ast(\tau^\prime,\tau,F)$ and consider their volume elements
$\alpha(\tau^\prime,\tau,E)\in\det(D_\ast(\tau^\prime,\tau,E))$ and
$\alpha(\tau^\prime,\tau,F)\in\det(D_\ast(\tau^\prime,\tau,F))$.
>From formula (5) for their determinants
it is clear that $\phi$  defines an isomorphism
$$\tilde\phi:\det(D_\ast(\tau^\prime,\tau,E))\to 
\det(D_\ast(\tau^\prime,\tau,F))$$
which (as above) is the product of the maps 
$$\Gamma(\Delta^\prime,\det(E))\to \Gamma(\Delta^\prime,\det(F))
\quad\text{and}\quad
\Gamma(\Delta,\det(E))\to \Gamma(\Delta,\det(F))$$ 
induced by $\phi$ or 
their contragradient maps. 

We claim now that {\it the constructed map $\tilde\phi$ preserves the volume;}
in other words,
$$\tilde\phi(\alpha(\tau^\prime,\tau,E))=\alpha(\tau^\prime,\tau,F).\tag6$$    
We will prove this fact assuming that $\tau^\prime$ is obtained
from $\tau$ by dividing one $q$-dimensional cell $e$ of
$\tau$ into two $q$-dimensional cells $e_+$ and $e_-$ 
and introducing a $(q-1)$-dimensional cell $e_0$, cf. figure 1. 
The general
statement follows from this special case by induction.

\midspace{4cm}\caption{Figure 1}

Under this 
our assumption the chain complexes $D_\ast(\tau^\prime,\tau,E)$ and
$D_\ast(\tau^\prime,\tau,F)$ are particularly simple. They have only
chains in dimensions $q-1$ and $q$ and can be described as follows.
The space $D_q=D_q(\tau^\prime,\tau,E)$ can be identified with
$D_q=\Gamma(e_+,E)$, the $q-1$-dimensional chains $D_{q-1}$ can be
identified with $D_{q-1}=\Gamma(e_0,E)$ and the boundary homomorphism
$\partial: D_q\to D_{q-1}$ can be identified with the restriction map
$r: \Gamma(e_+,E)\to \Gamma(e_0,E)$. Thus the volume
$$\alpha(\tau^\prime,\tau,E)\in\det(D_\ast(\tau^\prime,\tau,E))$$
is represented by the determinant of the restriction map $r$
$$\det(r):\Gamma(e_+,\det(E))\to \Gamma(e_0;\det(E)).$$
A similar description is valid for the bundle $F$. Now, the claim that
$\tilde\phi(\alpha(\tau^\prime,\tau,E))=\alpha(\tau^\prime,\tau,F)$
is precisely equivalent to the commutativity of the diagram
$$
\CD
\Gamma(e_+,\det(E))@>r>>\Gamma(e_0,\det(E))\\
@V\phi VV                @V\phi VV \\   
\Gamma(e_+,\det(F))@>r>>\Gamma(e_0,\det(F))
\endCD
$$
which is obvious. This finishes the proof that $\tilde \phi$ is volume
preserving.

Let us now define a map
$$S_{\tau^\prime,\tau}:\det(C_\ast(K,\tau^\prime,E))\to 
\det(C_\ast(K,\tau,E)).\tag7$$
We have the canonical identification 
$$U:\det(C_\ast(K,\tau^\prime,E))\to \det(C_\ast(K,\tau,E))\otimes
\det(D_\ast(K,\tau^\prime,\tau)).$$
For $x\in\det(C_\ast(K,\tau^\prime,E))$ set
$$S_{\tau^\prime,\tau}(x)=U(x)/\alpha(\tau^\prime,\tau,E).\tag8$$

Then, on one hand, 
we have the following commutative diagram
$$
\CD
\det(C_\ast(K,\tau^\prime,E))@>{S_{\tau^\prime,\tau}}>> \det(C_\ast(K,\tau,E))\\
        @V{T_{\tau^\prime}}VV                               @VV{T_\tau}V\\ 
\det(H_\ast(K,E))@>{=}>>                        \det(H_\ast(K,E)).
\endCD
\tag9
$$
To justify (9) one may refer to \cite{F}; theorem 1.3 of \cite{F} implies
(9) immediately.
On the other hand, for any bundle isomorphism $\phi:\det(E)\to\det(F)$ we
have the commutative diagram
$$
\CD
\det(C_\ast(K,\tau^\prime,E))@>{S_{\tau^\prime,\tau}}>> \det(C_\ast(K,\tau,E))\\ 
@V{\tilde\phi}VV                                        @V{\tilde\phi}VV\\
\det(C_\ast(K,\tau^\prime,F))@>{S_{\tau^\prime,\tau}}>> \det(C_\ast(K,\tau,F)). 
\endCD
\tag10
$$
It follows from the definition (8) and the fact established above (cf. (6)) 
that the induced map 
$\tilde\phi:\det(D_\ast(\tau^\prime,\tau,E))\to 
\det(D_\ast(\tau^\prime,\tau,F))$
preserves the volumes:
$\tilde\phi(\alpha(\tau^\prime,\tau,E))=\alpha(\tau^\prime,\tau,F).$    

This proves independence of the constructed map $\hat\phi$ on the 
triangulation. $\square$
\enddemo

\proclaim{1.3. Proposition (Functorial property)} 
Suppose that $E, F, \text{and}
\ G$ are three flat vector bundles over $K$ and let $\phi:\det(E)\to\det(F)$
and $\psi:\det(F)\to\det(G)$ be two isomorphisms. Then 
$$\widehat{\psi\circ\phi}=\hat \psi\circ\hat\phi:
\det(H_\ast(M;E))\to\det(H_\ast(M;G))\tag11$$
\endproclaim
\demo{Proof} It is clear from the construction.
\enddemo
\subheading{1.4. Remark} Suppose 
that $\phi:\det(E)\to\det(F)$ is an 
isomorphism
of flat vector bundles. Let $\lambda$ be a number.
Then
$$\widehat{\lambda\phi}=\lambda^{\chi(K)}\hat\phi,\tag12$$
where $\chi(K)$ denotes the Euler characteristic of $K$.

It follows from this remark that under the assumption that $\chi(K)=0$, 
the isomorphism $\hat\phi$ does not depend on $\phi$ (it only requires
existence of an isomorphism $\det(E)\to\det(F)$). 

\subheading{1.5} Correspondces between determinant lines appear also 
in the following version.

Let $E$ be a flat vector bundle (over the field $k\ =\ \C \ \text{or}\ \R$)
and let $\det_+(E)$ denote the flat bundle with fibre $\R_+$ obtained as 
follows. First consider the complement of the zero section of $\det(E)$
as a principal $k^\ast$-bundle; then form the associated bundle with fibre
$\R_+$ corresponding to the action of $k^\ast$ on $\R_+$ given by
$\alpha\cdot x\ =\ |\alpha|\cdot x$ for $\alpha\in k^\ast$, $x\in \R_+$.

Suppose now that we have two flat vector bundles $E$ and $F$ over a 
polyhedron $K$ and an isomorphism $\phi:\det_+(E)\to\det_+(F)$ of flat 
bundles. In terms of the flat bundles 
$\det(E)$ and $\det(F)$ this means that for simply-
connected open sets $U\subset K$ $\phi$ determines an isomorphism
$$\Gamma(U,\det(E))\to \Gamma(U,\det(F))$$
up to multiplication by a number with norm 1. Repeating the 
construction of section 1.1 we obtain that $\phi$ {\it determines a 
correspondence
$$\hat\phi : \det H_\ast(K;E)\to \det H_\ast(K;F)$$  
up to multiplication by a number with norm 1}.

Yet another version of correspondences between the 
determinant lines can be constructed as follows.
     
\proclaim{1.6. Proposition} Suppose that $E$ and $F$ are two flat vector 
bundles over a compact polyhedron $K$. Any flat section (metric) 
on the flat line bundle
$$\det(E)\otimes \det(F)$$
determines canonically an element of (or a metric on, correspondingly) 
the product of the determinant lines
$$\det(H_\ast(K,E))\otimes \det(H_\ast(K,F)).$$
\endproclaim
\demo{Proof} The construction is quite similar to the one described above.
We first remark that the line 
$\det(H_\ast(K,E))\otimes \det(H_\ast(K,E^\ast))$
can be identified with the product
$$\prod [\Gamma(\Delta,\det(E))\otimes 
\Gamma(\Delta, \det(F))]^{\epsilon (\Delta)}$$
and each term of this product has a section (or a metric) 
induced by the data. Then one
shows (using arguments similar to given above) that the constructed metric 
does not depend on the particular triangulation.
$\square$
\enddemo

\subheading{1.7. Remark} The metric constructed in Proposition 1.6 can 
also be 
obtained using the standard construction of Reidemeister metric 
(cf. 4.5) by first observing that $\det(E\oplus F)\ =\ \det(E)\otimes\det(F)$ 
and so the data determine a flat metric on this bundle; then the 
standard construction produces a metric on the determinant line
$\det(H_\ast(K,E\oplus F)$ which is canonically isomorphic to the product
$\det(H_\ast(K,E)\otimes \det(H_\ast(K,F))$. It is easy to see that the 
constructed metric coincides with the one given by Proposition 1.6.

\proclaim{1.8. Corollary} For any flat vector bundle $E$ over the 
polyhedron $K$
there is a canonical element in the line
$$\det(H_\ast(K,E))\otimes \det(H_\ast(K,E^\ast)),$$
(defined up to a sign) where $E^\ast$ denotes the dual flat vector bundle.
This element determines a canonical metric (which will be denoted
by $<\cdot>_E$) on the above line.
\endproclaim

The next proposition describes the relation between two constructions
which apeared in this section: between the correspondences $\hat\phi$ 
and the metrics $<\cdot>_E$
on the products of the determinant lines.

\proclaim{1.9. Proposition} Suppose that $E$ and $F$ are two flat vector 
bundles 
over $K$ and let $\phi:\det(E)\to\det(F)$ be an isomorphism of flat bundles.
Denote by $\psi:\det(F^\ast)\to\det(E^\ast)$ the adjoint bundle isomorphism.
Then for any
$$x\in \det(H_\ast(K,E))\qquad\text{and}\qquad y\in \det(H_\ast(K,F^\ast))$$
the following formula holds:
$$<x\otimes\hat\psi(y)>_E\quad =\quad <\hat\phi(x)\otimes y>_F;\tag13$$
here 
$\hat\phi : \det H_\ast(K;E)\to \det H_\ast(K;F)$ and
$\hat\psi : \det(H_\ast(K;F^\ast))\to \det(H_\ast(K;E^\ast))$ are constructed
as explained in 1.1 and $<\cdot>_E$ and $<\cdot>_F$ denote 
the canonical metrics
on the products $\det H_\ast(K;E)\otimes \det H_\ast(K;E^\ast)$ and
$\det H_\ast(K;F)\otimes \det H_\ast(K;F^\ast)$ respectively, constructed
in 1.8. $\square$
\endproclaim

\heading 2. Correspondences between cohomological determinants\endheading

2.1. Let $K$ be a finite polyhedron and let $E$ be a (real or complex)
flat vector bundle over
$K$. To simplify our notations the determinant of the homology of $K$ 
with coefficients in $E$ will be denoted
$$L_\bullet(E) = \det(H_\ast(K,E)).$$
Similarly, we will denote by
$$L^\bullet(E) = \det(H^\ast(K,E))$$   
the determinant of the cohomology; the latter can be understood as the 
cohomology of the locally constant sheaf determined by $E$.

\subheading{2.2} There is the canonical pairing
$$\{\ ,\ \}_E: L_\bullet(E)\otimes L^\bullet(E^\ast) \to \C\tag14$$ 
which is
determined by the evaluation maps
$H_i(K;E)\otimes H^i(K;E^\ast)\to \C$.
This gives a natural isomorphism
$$L^\bullet(E)\to L_\bullet(E^\ast)^\ast;\tag15$$
Using it we will reformulate the results of the previous section
using cohomological determinants.

\subheading{2.3} Suppose that $E$ and $F$ are two flat vector bundles and let
$\phi:\det(E)\to\det(F)$ be an isomorphism of flat bundles. Consider
the adjoint $\phi^\ast$ which acts  
$\phi^\ast: \det(F^\ast)\to\det(E^\ast)$; via the construction above,
it induces the correspondence
$$\widehat{(\phi^\ast)} : L_\bullet(F^\ast)\to L_\bullet(E^\ast)$$
and the adjoint of the latter gives the correspondence
$L_\bullet(E^\ast)^\ast\to L_\bullet(F^\ast)^\ast$. Using the canonical
isomorphism (14) we can interpret the last map as
$$\Check \phi : L^\bullet(E)\to L^\bullet(F).\tag16$$
We will refer to it as to the {\it cohomological correspondence determined by
$\phi$}.

\subheading{2.4} The cohomological correspondence has 
properties similar to the homological one.
Namely, 
\roster
\item {\it If $\chi(K)\ =\ 0$ then 
$\Check\phi: L^\bullet(E)\to L^\bullet(F)$
does not depend on $\phi$.
\item suppose that $E, F, \text{and}\ G$ are three flat vector bundles 
over $K$ and let $\phi:\det(E)\to\det(F)$
and $\psi:\det(F)\to\det(G)$ be two isomorphisms. Then 
$$\Check {\psi\circ\phi}=\Check \psi \circ  \Check \phi: L^\bullet(E)\to L(G^\bullet)$$

\item any flat section (metric) on the flat line bundle
$\det(E)\otimes\det(F)$ determines canonically an element (or a metric,
correspondingly) on the product
$L^\bullet(E)\otimes L^\bullet(F).$}  
\newline
Indeed, any flat section (metric) on the line bundle 
$\det(E)\otimes\det(F)$ determines 
the obvious metric on the dual bundle
$\det(E^\ast)\otimes\det(F^\ast)$ and the latter by the construction of
proposition 1.6
gives an element (metric, correspondingly) on the line
$L_\bullet(E^\ast)\otimes L_\bullet (F^\ast)$. The latter
clearly determines an element
(or metric) on $L^\bullet(E)\otimes L^\bullet(F).$   
\newline
In particular, {\it for any flat line bundle $E$ we obtain a canonical metric 
on the line $L^\bullet(E)\otimes L^\bullet(E^\ast)$; this metric will be
denoted} $<\cdot>_E$.

\item[4] There is also the relation similar to (13).
{\it Suppose that $E$ and $F$ are two flat vector bundles 
over $K$ and let $\phi:\det(E)\to\det(F)$ be an isomorphism of flat bundles.
Denote by $\psi:\det(F^\ast)\to\det(E^\ast)$ the adjoint bundle isomorphism.
Then we have the combinatorial correspondences
$\Check \phi: L^\bullet(E)\to L^\bullet(F)$ and 
$\Check \psi: L^\bullet(F^\ast)\to L^\bullet(E^\ast),$ and for any
$x\in L^\bullet(E)$ and $y\in L^\bullet(F^\ast)$ the following formula
holds}
$$<x\otimes\Check\psi(y)>_E\quad =\quad <\Check\phi(x)\otimes y>_F.\tag17$$
\endroster

\heading 3. Poincar\'e duality for determinant lines\endheading

In this section $K$ {\it will denote a 
closed oriented piecewise-linear manifold of odd dimension}. $E$ will
denote a flat $k$-vector bundle over $K$, where $k$ is $\R$ or $\C$.

It is well known that the Poincar\'e duality induces some isomorphisms
between the determinant lines, cf. \cite{M1}, \cite{W}, \cite{BZ}.
Our aim in this section is to introduce the notations and
establish relations between the correspondences
$\hat\phi$ and $\Check\phi$ introduced in the previous sections and the 
isomorphisms determined by the Poincar\'e duality. The result of this section
will be used in \S 4 to define the Poincar\'e-Reidemeister metric.

\subheading{3.1} Since the orientation of the manifold $K$ is 
supposed to be fixed, for any integer $q$ we have a nondegenerate 
intersection form
$$H_q(K;E)\otimes H_{n-q}(K,E^\ast)\to k,$$
where $n=\dim K$ is supposed to be odd. 
The above pairing allows one to identify
$H_q(K,E)$ with the dual of $H_{n-q}(K,E^\ast)$. Since $q$ and $n-q$ are
of opposite parity, it defines an isomorphism 
$$\D_E: \ L_\bullet(E)\to L_\bullet(E^\ast).\tag18$$

\subheading{3.2} Another description of the Poincar\'e duality 
map $\D_E$ can 
be given as follows. Consider a triangulation 
$\tau$ of $K$ and
the dual cell decomposition $\tau^\ast$. There is the intersection
pairing on the chain level
$$C_q(K,\tau,E)\otimes C_{n-q}(K,\tau^\ast,E^\ast)\to k.$$
Let $\Delta$ be a $q$-dimensional simplex of $\tau$ and let $\Delta^\ast$
be the dual $(n-q)$-dimensional cell of $\tau^\ast$. 
Then the above 
intersection form splits into an orthogonal sum of partial pairings
$$\Gamma(\Delta,E)\otimes \Gamma(\Delta^\ast,E^\ast)\to k,$$
which assign to a pair of flat sections $s\in\Gamma(\Delta, E)$ and
$s^\prime\in\Gamma(\Delta^\ast, E^\ast)$ the number 
$<s^\prime(x), s(x)>\in k$,
where $x$ denotes the common point of $\Delta$ and $\Delta^\ast$.
The last pairing is non-degenerate and defines an isomorphism
$$\alpha_\Delta:\Gamma(\Delta,E)\to\Gamma(\Delta^\ast,E^\ast)^\ast.$$
Thus we obtain the maps
$$\det(\alpha_\Delta): \Gamma(\Delta,\det(E))\to
\Gamma(\Delta^\ast,\det(E^\ast))^\ast.\tag19$$    

We claim that the following diagram commutes
$$
\CD
\prod_\Delta
\Gamma(\Delta,\det(E))^{\epsilon(\Delta)}
@>{\otimes\det(\alpha_\Delta)^{\epsilon(\Delta)}}>>
\prod_{\Delta^\ast}\Gamma(\Delta^\ast,\det(E^\ast))^{\epsilon(\Delta^\ast)}\\ 
@V{T_\tau}VV              @VV{T_{\tau^\ast}}V\\
L_\bullet(E)@>{\D_E}>>  L_\bullet(E^\ast)\ .
\endCD
\tag20$$
This gives another characterization of the duality map
$\D_E: L_\bullet(E)\to L_\bullet(E^\ast)$. 

The claim is a special case of the following general algebraic remark. 
Let $C$ be a chain complex of finite dimensional vector spaces and let
$n$ be an odd integer. Form a new complex $D\ =\ (D_j,\delta)$, where
$D_j\ =\ C_{n-j}^\ast$ and $\delta: D_j\to D_{j-1}$ is defined to be
the dual of $\partial: C_{n-j+1}\to C_{n-j}$. For any index $i$
we have the canonical map $\alpha_i: C_i\to D^\ast_{n-i}$ (the identity)
and $\det(\alpha_i): \det(C_i)\to \det(D_{n-i})^\ast$. We obtain the 
following diagram
$$
\CD
\det(C)=\otimes\det(C_i)^{(-1)^i}@>{\otimes \det(\alpha_i)^{(-1)^i}}>>
\otimes\det(D_i)^{(-1)^i}=\det(D)\\
@V{T_C}VV         @V{T_D}VV\\
\det(H_\ast(C))@>\Sigma>>\det(H_\ast(D))
\endCD
\tag21$$
where the vertical maps $T_C$ and $T_D$ are the canonical maps and  the 
horizontal map $\Sigma$ is determined by the obvious maps
$H_i(C)\to H_{n-i}(D)^\ast$. The proof of this elementary fact can be
obtained directly from the definitions.


\proclaim{3.3. Proposition} The duality isomorphism $\D_E$
has the following properties:
\roster
\item  The map 
$$\D_{E^\ast}: L_\bullet(E^\ast)\to L_\bullet(E^{\ast\ast})=L(E)$$
is inverse to $D_E$;
\item let $E$ and $F$ be two flat vector bundles over $K$ and let
$\phi:\det(E)\to \det(F)$ be an isomorphism of flat bundles. Denote by
$\psi:\det(F)^\ast=\det(F^\ast)\to \det(E^\ast)=\det(E)^\ast$
the map adjoint to $\phi$. Then there are two combinatorial correspondences
$$\hat\phi:L_\bullet(E)\to L_\bullet(F)\quad\text{and}\quad
\hat\psi:L_\bullet(F^\ast)\to L_\bullet(E^\ast)$$
and the following diagram
$$
\CD
L_\bullet(E)@>{\hat\phi}>>L_\bullet(F)\\
@V{\D_E}VV @V{\D_F}VV\\
L_\bullet(E^\ast)@<{\hat{\psi}}<<L_\bullet(F^\ast)
\endCD
\tag22$$
commutes.
\endroster
\endproclaim
\demo{Proof} The first property above follows immediately from the
symmetry of the intersection numbers.

The second property follows from the second description of the map $\D_E$
given above
and from the commutative diagram
$$
\CD
\Gamma(\Delta,\det(E))@>{\det(\alpha_\Delta)}>>
(\Gamma(\Delta^\ast,\det(E^\ast)))^\ast\\
@V{\phi}VV                        @V{\psi^\ast}VV\\
\Gamma(\Delta,\det(F))@>{\det(\alpha_\Delta)}>>
(\Gamma(\Delta^\ast,\det(F^\ast)))^\ast. 
\endCD
$$
This completes the proof. $\square$
\enddemo

\subheading{3.4} For any closed oriented manifold $K$ of odd dimension
and a flat vector bundle $E$ over $K$ 
there is the Poincar\'e duality map
$$\D_E: L^\bullet(E) \to L^\bullet(E^\ast)$$
acting on cohomological determinant lines. It can be defined as follows.
Using the canonical pairings (13)
$$\{\ ,\ \}_E: L_\bullet(E)\otimes L^\bullet(E^\ast)\to \C$$
and the similar pairing for the dual flat vector bundle $E^\ast$,
we will define the map
$\D_{E^\ast}: L^\bullet(E^\ast) \to L^\bullet(E^{\ast\ast})= L^\bullet(E)$  
by the requirement
$$\{\D_E(x)\otimes \D_{E^\ast}(y)\}_{E^\ast}\ =\ \{x\otimes y\}_E\tag23$$
for any $x\in L_\bullet(E)$ and $y\in L^\bullet(E^\ast)$.

It is clear that the above 
duality map $\D_E: L^\bullet(E) \to L^\bullet(E^\ast)$
can also be described as the map between the determinant lines of the 
cohomology induced by the isomorphisms
$$H^i(K,E)\ \to\ (H^{n-i}(K,E^\ast))^\ast,$$
comming from the nondegenerate intersection forms
$H^i(K,E)\otimes H^{n-i}(K,E^\ast)\to k.$
\subheading{3.5} Note that if the dimension of the manifold $K$ is even 
then the Poincar\'e duality determines an element of the line
$$\D_E\in L^\bullet(E)\otimes L^\bullet(E^\ast)$$
It can be shown that 
$$<\D_E>_E\ =\ 1;$$
thus, in the even-dimensional case the Poincar\'e duality $\D_E$
determines the canonical pairing $<\ >_E$.
Compare \cite{W}.

\heading 4. Poincar\'e - Reidemeister metric\endheading

In this section $K$ will denote 
an {\it odd-dimensional closed oriented manifold}. We will construct  
canonical metrics, which we call Poincar\'e - Reidemeister metrics, 
on the determinant lines $L_\bullet(E)$ and $L^\bullet(E)$ for any flat
vector bundle.
\subheading{4.1} Let
$E$ be a flat vector bundle over $K$. By Corollary 1.8
there is a canonical metric on the line 
$L_\bullet(E)\otimes L_\bullet(E^\ast)$ which we will
denote by $<\cdot>_E$. On the other hand, there is defined the duality
map
$$\D_E: L_\bullet(E)\to L_\bullet(E^\ast),$$
cf. 3.1. 
\subheading{Definition} The Poincar\'e - Reidemeister metric 
$||\cdot||_{L_\bullet(E)}^{PR}$ on the line $L_\bullet(E)$ is defined by
$$||x||_{L_\bullet(E)}^{PR}\ =\ <x\otimes \D_E(x)>^{1/2}_E\tag24$$
for $x\in L(E)$.

First, we want to show that the Poincar\'e - Reidemeister metric behaves well
with respect to the correspondences $\hat\phi$ between the determinant lines.

\proclaim{4.2. Proposition} Suppose that $K$ is a closed oriented manifold 
of
odd dimension. Let $E$ and $F$ be two flat bundles over $K$ and let
$\phi:\det(E)\to\det(F)$ be an isomorphism of flat bundles. 
Let $\hat\phi:L_\bullet(E)\to L_\bullet(F)$ be the induced correspondence
between the determinant lines, cf. 1.1. Then the correspondence $\hat\phi$
preserves the Poincar\'e - Reidemeister metrics:
for any
$x\in L_\bullet(E)$ we have
$$||\hat\phi(x)||_{L_\bullet(F)}^{PR}\ =\ ||x||_{L_\bullet(E)}^{PR}\tag25$$
where $||\cdot||_{L_\bullet(E)}^{PR}$ and $||\cdot||_{L_\bullet(F)}^{PR}$ 
denote the Poincar\'e - Reidemeister
metrics on $L_\bullet(E)$ and $L_\bullet(F)$ correspondingly.
\endproclaim
\demo{Proof} According to the definitions and using Propositions 1.9 and
3.3.(2) we obtain
$$
\CD
||\hat\phi(x)||_{L_\bullet(F)}^{PR,2}\ =\ <\hat\phi(x)\otimes\D_F\hat\phi(x)>_F\ =\\
<x\otimes \hat\psi\D_F\hat\phi(x)>_E\ =\\
<x\otimes \D_E(x)>_E\  =\ ||x||^{PR,2}_{L_\bullet(E)}
\endCD
$$
where $<\cdot>_E$ and $<\cdot>_F$ denote the canonical metrics on the lines
$L_\bullet(E)\otimes L_\bullet(E^\ast)$ and 
$L_\bullet(F)\otimes L_\bullet(F^\ast)$ respectively.
\enddemo

\subheading {4.3} Now I will describe a recipe computing the Poincar\'e -
Reidemeister metric from the combinatorial data. 

Let $K$ denote a closed oriented odd-dimensional manifold and let $E$ be
a flat vector bundle over $K$. Let $\tau$ be a triangulation of $K$.
Given an element $\alpha\in \det C_\ast(K,\tau,E)$ our task is to compute
the Poincar\'e - Reidemeister metric of the corresponding element
$||T_\tau(\alpha)||^{PR}_{L_\bullet(E)}$, where
$T_\tau:\det C_\ast(K,\tau,E)\to \det H_\ast(K;E)$ is the canonical 
isomorphism. 

Let $\tau^\ast$ be the dual cell decomposition of $K$ and let $E^\ast$
be the dual flat vector bundle. Choose two arbitrary non-zero elements
$$\beta\in \ \det C_\ast(K,\tau^\ast,E)$$
and
$$\gamma \in\ \det C_\ast(K,\tau^\ast,E^\ast).$$
We are going to define three positive real numbers
$$[\alpha]/[\beta],\quad <\beta,\gamma>,\quad\text{and}\quad \alpha/\gamma.$$
The first number $[\alpha]/[\beta]=\lambda >0$ is defined by the requirement
$$T_\tau(\alpha)\ =\ \pm \lambda T_{\tau^\ast}(\beta),$$
the equality taking place in the determinant line $\det H_\ast(K;E)$.

The second number is defined by
$$ <\beta,\gamma> \ =
\ \prod_{\Delta^\ast}\{\beta_{\Delta^\ast}, 
\gamma_{\Delta^\ast}\}_{\Delta^\ast},$$
where the product is taken over all cells $\Delta^\ast$ 
of the dual cell decosition $\tau^\ast$ and where
$$\beta\ =\ \otimes\beta_{\Delta^\ast} \ \in\ 
\otimes \Gamma(\Delta^\ast, \det(E))^{\epsilon(\Delta^\ast)}$$
and
$$\gamma\ =\ \otimes\gamma_{\Delta^\ast} \ \in\ 
\otimes \Gamma(\Delta^\ast, \det(E^\ast))^{\epsilon(\Delta^\ast)}.$$
Here the brackets $\{\ ,\ \}_{\Delta^\ast}$ denote the absolute value of the
canonical evaluation pairing on
$$\Gamma(\Delta^\ast, \det(E))^{\epsilon(\Delta^\ast)}\otimes
\Gamma(\Delta^\ast, \det(E^\ast))^{\epsilon(\Delta^\ast)}.$$
                                                                
The third number $\alpha/\gamma$ is defined by the Poincar\'e duality: if
$\alpha=\otimes \alpha_\Delta\in \Gamma(\Delta,\det(E))^{\epsilon(\Delta)}$
then
$$ \alpha/\gamma \ =\ \prod_{(\Delta,\Delta^\ast)}
(\alpha_\Delta: \gamma_{\Delta^\ast})$$
where the product is taken over all pairs $(\Delta,\Delta^\ast)$ of mutually
dual cells and for any such pair the number 
$(\alpha_\Delta: \gamma_{\Delta^\ast})$ is the absolute value of the 
evaluation at the common point of $\Delta$ and $\Delta^\ast$, cf. figure 2,
of the flat sections 
$\alpha_\Delta\in \Gamma(\Delta,\det(E)^{\epsilon(\Delta)}$ 
and 
$\beta_{\Delta^\ast}\in \Gamma(\Delta,\det(E^\ast)^{\epsilon(\Delta^\ast)+1}.$ 

\midspace{4cm}\caption{Figure 2}

Now we can define the Poincar\'e - Reidemeister metric as the product of the 
constructed three numbers:
$$
||T_\tau(\alpha)||^{PR,2}_{L_\bullet(E)}\ =\ 
[\alpha]/[\beta]\ \times \ <\beta,\gamma>\ \times\ \alpha/\gamma\tag26$$
The above definition clearly does not depend on the choice of $\beta$ and 
$\gamma$ and the result coincides with the definition (24).

\proclaim{4.4. Proposition} Suppose that flat vector bundle $E$ is 
unimodular,
i.e. there exists a flat metric on the bundle $\det(E)$. Then 
there is defined the standard Reidemeister metric on
the determinant line $L_\bullet(E)$. 
We claim that in this case this Reidemeister metric 
coincides with the Poincar\'e - Reidemeister metric constructed above.
\endproclaim

\subheading{4.5} Let us first recall the construction of 
the Reidemeister 
metric assuming that $\det(E)$ is unimodular.
Let $\mu$ be a flat metric on the 
flat bundle $\det(E)$. Fix a polyhedral cell decomposition 
$\tau$ of $K$.
For any cell $\Delta\subset K$ the flat metric $\mu$ defines an element
$$\mu_\Delta\in\det \Gamma(\Delta,E)=\Gamma(\Delta,\det(E)),$$
determined up to multiplication by a number with norm 1.
The product of these $\mu_\Delta$'s 
$$\mu_{\tau,E}=\prod_\Delta \mu_\Delta^{\epsilon(\Delta)},\qquad 
\epsilon(\Delta)=(-1)^{\dim(\Delta)}$$
defines an element  
$$\mu_{\tau,E}\in\det C_\ast(K,\tau,E)=
\prod\det\Gamma(\Delta,E)^{\epsilon(\Delta)}$$
and thus it defines an element $T_\tau(\mu_{\tau,E})\in L_\bullet(E)$ 
of the determinant line, determined again up to multiplication by a number
with norm 1.
The later element correctly defines a metric 
$||\cdot||_{L_\bullet(E)}$ on the determinant line $L_\bullet(E)$ 
by the requirement 
that $||T_\tau(\mu_{\tau,E})||_{L_\bullet(E)}=1$. 
If one choses another flat metric $\mu^\prime =\lambda\mu$ on $\det(E)$
then the constructed metric on $L_\bullet(E)$ does not change, since
$\chi(K)=0$ by the Poincar\'e duality.

One of the main properties of the Reidemeister metric is its 
{\it combinatorial invariance}. This means that for a pair of
polyheral cell decompositions 
$\tau$ and $\tau^\prime$ of $K$ having a common subdivision
we have the diagram
$$
\CD
\mu_{\tau,E}\in\det(C_\ast(K,\tau,E))@>{T_\tau}>>
\det H_\ast(K,E)@<{T_{\tau^\prime}}<<
\det(C_\ast(K,\tau^\prime,E))
\ni\mu_{\tau^\prime,E}
\endCD
$$
and
$$T_\tau(\mu_{\tau,E})\ :\ T_{\tau^\prime}(\mu_{\tau^\prime,E})$$
is a number of norm 1.

\subheading{4.6. Proof of Proposition 4.4} Let $\mu$
and $\nu$ be a pair of mutually dual flat metrics on $\det(E)$ and 
$\det(E^\ast)$ correspondingly.
Let $\tau$ be a triangulation of $K$ and let $\tau^\ast$ be the 
dual cell decomposition.
Let $\mu_{\tau,E}\in L_\bullet(E)$ be the corresponding element 
as defined above in 4.5. 
To prove our statement we have to show that
$$<T_\tau(\mu_{\tau,E})\otimes \D_E T_\tau(\mu_{\tau,E})>_E\ =1.$$

First note that 
$$T_\tau(\mu_{\tau,E})=T_{\tau^\ast}(\mu_{{\tau^\ast},E})$$
by the combinatorial invariance property of the Reidemeister metric under 
subdivisions, cf. 4.5. On the other hand
$$\D_E T_{\tau^\ast}(\mu_{{\tau^\ast},E})=T_\tau(\nu_{\tau,E^\ast}),$$
as one directly checks using the definitions.
Thus we have
$$
\CD
<T_\tau(\mu_{\tau,E})\otimes \D_E T_\tau(\mu_{\tau,E})>_E\ =\\
<T_\tau(\mu_{\tau,E})\otimes T_\tau(\nu_{\tau,E^\ast})>_E\ =1
\endCD
$$   
completing the proof.
$\square$

\subheading{4.7} Let us now define a similar 
{\it Poincar\'e - Reidemeister metric} 
(which we will denote $||\cdot||_{L^\bullet(E)}^{PR}$) on the cohomological
determinant $L^\bullet(E)$. Here again $K$ denotes
a closed oriented odd-dimensional manifold and $E$ a flat vector bundle
over $K$. 

Recall that there is a canonical norm $<\ ,\ >_E$ on the line
$L^\bullet(E)\otimes L^\bullet(E^\ast)$  (cf. 2.4.(3)) and there is the 
Poincare duality map $\D_E: L^\bullet(E)\to L^\bullet(E^\ast)$
(cf. 3.4). 
\subheading{Definition} For $x\in L^\bullet(E)$ define
$$||x||_{L^\bullet(E)}^{{PR},2}\ =\ <x\otimes \D_E(x)>\tag27$$

Note that this Poincar\'e - Reidemeister metric on the cohomological 
determinant
line is again a purely {\it combinatorial} invariant of $K$ and $E$.

The cohomological version of the Poincar\'e - Reidemeister metric has 
properties
similar to the homological one. We will formulate them without proof:

\proclaim{4.8. Proposition}\roster
\item If a flat bundle $E$ is unimodular (i.e. the flat line bundle
$\det(E)$ admits a flat metric) then the Poincar\'e - Reidemeister
metric on $L^\bullet(E)$ coincides with the standard Reidemeister metric;
\item Suppose that $E$ and $F$ are two flat bundles over $K$ and let
$\phi:\det(E)\to\det(F)$ be an isomorphism of flat bundles. Let
$\Check \phi:L^\bullet(E)\to L^\bullet(F)$ 
denote the induced correspondence
between the determinant lines, cf.2.3. Then the correspondence $\Check\phi$
preserves the Poincar\'e - Reidemeister metrics; in other words, for
any $x\in L^\bullet(E)$ we have
$$||\Check\phi(x)||_{L^\bullet(F)}^{PR}\ =\ 
||x||_{L^\bullet(E)}^{PR}.\tag28$$
\endroster
\endproclaim  

\heading 5. Ray-Singer norm of the combinatorial
correspondence between the determinant
lines\endheading

Let $K$ denote a {\it closed oriented smooth manifold}.

Suppose that we are given two flat vector bundles $E$ and $F$ over $K$
and an isomorphism $\phi:\det(E)\to\det(F)$ between their flat determinant
line bundles. As we have seen in section 1, there is
the correspondence
$$\Check\phi:L^\bullet (E)\to L^\bullet (F)$$
which is determined completely by the combinatorial structure of the
initial data $(K,E,F,\phi)$. The correspondence
$\Check\phi$ determines a metric on the relative determinant line
$L^\bullet (F):L^\bullet (E)$, which is defined by the requirement 
that the norm of $\Check\phi$ is 1. 

In this section we use the main theorem of J.-M.Bismut and W.Zhang 
\cite{BZ},
theorem 0.2, to show that the Ray-Singer construction of analytic torsion
produces the same metric on this relative determinant. 

In the next section we will show that this theorem about the metrics on
the relative determinants gives in the 
odd-dimensional situation an identification of the Ray-Singer metric with 
a combinatorially defined Poincar\'e - Reidemeister metric, defined in \S 4.

\proclaim{5.1. Theorem} Let $K$ be a closed oriented smooth manifold and let
$E$ and $F$ be two flat vector bundles over $K$ supplied with an
isomorphism of flat bundles $\phi:\det(E)\to\det(F)$. Consider the
correspondence $\Check\phi:L^\bullet (E)\to L^\bullet (F)$ (cf. \S 1)
constructed by means of a smooth polyhedral cell decomposition of $K$. 
Fix a Riemannian metric on $K$. Fix Hermitian metrics 
on the flat vector bundles $E$ and $F$ in 
such a way that the induced metrics on $\det(E)$ and $\det(F)$ are
isomorphic via $\phi: \det(E)\to \det(F)$. 
There are then defined the Ray-Singer metrics
$||\cdot||^{RS}_{L^\bullet(E)}$ and $||\cdot||^{RS}_{L^\bullet(F)}$ 
on the determinant lines $L^\bullet(E)$ and $L^\bullet(F)$ respectively.
We claim that the 
relative Ray-Singer metric on the line
$\Hom(L^\bullet(E),L^\bullet(F)) = L^\bullet(F)\otimes L^\bullet(E)^\ast$ 
coincides with the metric, determined on this line 
by the correspondence $\Check\phi$. 
In other words, the correpondence $\Check\phi$
preserves the Ray-Singer metrics: for any
$x\in L^\bullet(E)$ the following formula holds
$$||\Check\phi(x)||^{RS}_{L^\bullet(F)}\ =\ ||x||^{RS}_{L^\bullet(E)}.\tag29$$
\endproclaim
\demo{Proof} The proof easily follows from theorem (0.2) of Bismut and
Zhang \cite{BZ}. If we would have known that the cell decomposition 
associated
with an arbitrary Morse function is a smooth polyhedral cell decomposition
(in the sense of \cite{RS1})
with respect to the smooth structure of $K$ (which is not true in
general) then we would be able to 
apply the theorem of Bismut and Zhang directly.
Instead, we will proceed as follows.

Fix a smooth triangulation of $K$ (cf. \cite{Mu2})
and consider its second
derived subdivision (cf. [RS1]) which we will denote $\tau$. 

Choose a Riemanian metric on the manifold $K$.
There is a Morse function $f$ on $K$ having the following properties:
\roster
\item Any open simplex $\Delta$ of the triangulation $\tau$ contains a unique
critical point $p_\Delta$ of $f$ having index $\dim(\Delta)$;
\item The unstable manifold of the critical point $p_\Delta$ coincides
with $\Delta$.
\endroster
Such function $f$ can be constructed by considering the handle decomposition
associated with the triangulation $\tau$ (cf. [RS1], Prop. 6.9), then by
smoothing the corners of the handles (similarly to [M3]) and then by 
using the standard correspondence between glueing handles, elementary
cobordisms and the Morse functions having precisely one critical point,
cf. [M3].

The Thom-Smale complex associated with the function $f$ 
is now identical to the simplicial chain complex of $K$ with respect the 
triangulation $\tau$. 

For any simplex $\Delta$ of $\tau$ fix a volume element in
$\det(\Gamma(\Delta,E))\ =\ \Gamma(\Delta,\det(E))$. Then fix the 
volume element in $\det(\Gamma(\Delta,F))\ =\ \Gamma(\Delta,\det(F))$
corresponding to the choice made for $E$ under the bundle isomorphism
$\phi:\det(E)\to\det(F)$. The above choices determine the metrics
$||\cdot||^{M,X}_{L^\bullet(E)}$ and $||\cdot||^{M,X}_{L^\bullet(F)}$
via the canonical isomorphisms
$$L^\bullet(E)\to \det(C_\ast(K,\tau,E^\ast))^{-1}\tag30$$ 
(cf. (15))
and similarly for the bundle $F$; 
these metrics are called Milnor metrics in \cite{BZ}.
According to our general construction in \S 1, with the coherent choices
(with respect to $\phi$) as above
made for the volume elements, we have
$$||x||^{M,X}_{L^\bullet(E)}\ =
\ ||\Check\phi(x)||^{M,X}_{L^\bullet(F)}\tag31$$
for $x\in L^\bullet(E)$.

Let us apply theorem (0.2) of Bismut and Zhang \cite{BZ} twice: once for
the metrics 
$||\cdot||^{RS}_{L^\bullet(E)}$ and 
$||\cdot||^{M,X}_{L^\bullet(E)}$ on 
${L^\bullet(E)}$ and the second time for the metrics
$||\cdot||^{RS}_{L^\bullet(F)}$ and
$||\cdot||^{M,X}_{L^\bullet(F)}$
on ${L^\bullet(F)}$. The right hand sides of the formula (0.8) of \cite{BZ}
in both cases will be the same because it depends only on the metric on
the manifold $K$ and the flat connection and the metric on the determinants
$\det(E)$ and $\det(F)$. Thus, subtracting, we obtain
$$
\frac{||x||^{RS}_{L^\bullet(E)}}{||x||^{M,X}_{L^\bullet(E)}}\ =\ 
\frac{||y||^{RS}_{L^\bullet(F)}}{||y||^{M,X}_{L^\bullet(F)}} \tag32
$$
for any $x\in L^\bullet(E)$ and $y\in L^\bullet(F)$. Now, set
$y\ =\ \Check\phi(x)$; then from (30) we obtain
$||\Check\phi(x)||^{RS}_{L^\bullet(F)}\ =\ ||x||^{RS}_{L^\bullet(E)}$ 
as stated.
$\square$
\enddemo

The following formulation is equivalent to theorem 5.1.

\proclaim{5.2. Theorem} Let $K$ be a closed oriented smooth manifold and let
$E$ and $F$ be two flat vector bundles over $K$ such that the 
flat bundle $\det(E)\otimes\det(F)$ admits a flat metric $\mu$.
Then by 2.4.(3) $\mu$ determines 
(combinatorially) a
metric $||\cdot||_\mu$ on the line $L^\bullet(E)\otimes L^\bullet(F)$.
Fix a Riemannian metric on $K$. Fix Hermitian metrics 
on the flat vector bundles $E$ and $F$ in 
such a way that the induced metrics on $\det(E)\otimes\det(F)$ is $\mu$.
There are then defined the Ray-Singer metrics
$||\cdot||^{RS}_{L^\bullet(E)}$ and $||\cdot||^{RS}_{L^\bullet(F)}$ 
on the determinant lines $L^\bullet(E)$ and $L^\bullet(F)$ respectively.
We claim that the product of these Ray-Singer metrics coincides with
the combinatorial metric $||\cdot||_\mu$; more precisely
for any $x\in L^\bullet(E)$ and $y\in L^\bullet(F)$ we have
$$||x||^{RS}_{L^\bullet(E)}\cdot||y||^{RS}_{L^\bullet(F)}\ =
\ ||x\otimes y||_\mu.$$   
\endproclaim

We will formulate now a very special case of the above theorem, which will
be used below. We will see later that this simple statement allows to 
identify completely the Ray-Singer metric in the odd-dimensional case.

\proclaim{5.3. Theorem} Let $K$ be a closed oriented manifold and let $E$ be
a flat vector bundle over $K$. By 2.4.(3) there is a canonical 
metric $<\cdot>_E$ on the line
$L^\bullet(E)\otimes L^\bullet(E^\ast)$ which is determined by the
combinatorial structure of $K$ and $E$. Choose a Riemannian metrix on $K$
and a Hermitian metric on the bundle $E$ (which is not supposed to be flat).
The latter determines a metric on the bundle $E^\ast$. The construction
of Ray and Singer produces now the metrics 
$||\cdot||^{RS}_{L^\bullet(E)}$ and $||\cdot||^{RS}_{L^\bullet(E^\ast)}$ 
on the lines $L^\bullet(E)$ and $L^\bullet(E^\ast)$ correspondingly.
Then their product is equal to the canonical combinatorial metric
$<\cdot>_E$ on $L^\bullet(E)\otimes L^\bullet(E^\ast)$. 
\endproclaim

\heading 6. Ray-Singer metric coincides with the Poincar\'e - Reidemeister
metric\endheading

We will prove in this section that for any closed oriented odd-dimensional
manifold $K$ and for any flat vector bundle $E$ over $K$ the Ray-Singer
norm on the determinant line of the cohomology $L^\bullet(E)$
(constructed by using any Riemannian metric on $K$ and a metric on $E$)
coincides with the combinatorially defined Poincar\'e - Reidemeister
metric.

We will start with the following lemma.

\proclaim{6.1. Lemma} Let $K$ be an odd-dimensional manifold and let $E$
be a flat vector bundle over $K$. Consider the map
$\D_E: L^\bullet(E)\to L^\bullet(E^\ast)$ 
determined by the Poincare duality cf. 3.4. 
Fix a Riemannian metric on $K$ and a
metric on the bundle $E$; the latter determines a metric on $E^\ast$.
Consider now the Ray-Singer metrics 
$||\cdot||^{RS}_{L^\bullet(E)}$ and
$||\cdot||^{RS}_{L^\bullet(E^\ast)}$
on the lines $L^\bullet(E)$ and $L^\bullet(E^\ast)$
determined by the choices made above. Then the map $\D_E$ preserves
the Ray-Singer metrics: for any $x\in L^\bullet(E)$ we have
$$||x||^{RS}_{L^\bullet(E)}\ =\ ||\D_E(x)||^{RS}_{L^\bullet(E^\ast)}$$
\endproclaim
\demo{Proof} This is well-known, cf. \cite{BZ}, page 35.
\enddemo

\proclaim{6.2. Theorem} Let $K$ be a closed oriented smooth
odd-dimensional manifold and let $E$
be a flat vector bundle over $K$. Consider the  
Poincar\'e - Reidemeister metric $||\cdot||^{PR}_{L^\bullet(E)}$  
on the determinant line of the 
cohomology $L^\bullet(E)$. Choose an arbitrary Riemannian metric on $K$  
and an arbitrary Hermitian metric on $E$ and consider the metric 
$||\cdot||^{RS}_{L^\bullet(E)}$ given by the construction of Ray and 
Singer \cite{RS}. Then these two metrics coincide.
\endproclaim
\demo{Proof} For $x\in L^\bullet(E)$ we obtain
$$
\align
||x||^{{PR},2}_{L^\bullet(E)} \ &= 
<x\otimes\D_E(x)>_E\ \quad\qquad\quad\text{(by definition (27))}\\
&=\ ||x||^{RS}_{L^\bullet(E)}\times ||\D_E(x)||^{RS}_{L^\bullet(E)}\quad  
\text{(by Theorem 5.3)}\\
&=\ ||x||^{RS}_{L^\bullet(E)}\times ||x||^{RS}_{L^\bullet(E)} \qquad\quad
\text{(by Lemma 6.1)}
\endalign
$$
This completes the proof $\square$.
\enddemo
\subheading{6.3. Remark} Theorem 6.2 together with Proposition 4.8.(1)
clearly generalize theorem of W.M\"uller \cite{Mu1}.
On the other hand, proof of Theorem 6.2 can be based on the theorem
of W.M\"uller \cite{Mu1} instead of using theorem of J.-M.Bismut and
W.Zhang \cite{BZ}. Namely, if $K$ is odd-dimensional then Theorem
5.3 follows from theorem of W.M\"uller since the bundle $E\oplus E^\ast$
is unimodular and the Reidemeister metric on $L^\bullet(E)\otimes
L^\bullet(E^\ast)$ coincides with the canonical metric $< \cdot>_E$. 
Using Proposition 4.8.(2) one can obtain also Theorems 5.1 and Theorem 5.2
in the odd-dimensional case.

\end